\input amstex
\magnification=\magstep1 
\baselineskip=13pt
\documentstyle{amsppt}
\vsize=8.7truein
\CenteredTagsOnSplits \NoRunningHeads
\def\vl{\operatorname{vol}}
\def\PP{\Cal P}
\def\AA{\Cal A}
\def\VV{\Cal V}
\def\EE{\bold{E\thinspace}}

\def\TT{\Cal T}

\def\xx{\bold{x}}
\def\yy{\bold{y}}
\def\sss{\bold{s}}
\def\ttt{\bold{t}}

\def\Pr{\bold{Pr\thinspace }}

\topmatter
\title What does a random contingency table
look like? \endtitle
\author Alexander Barvinok \endauthor
\address Department of Mathematics, University of Michigan, Ann Arbor,
MI 48109-1043, USA \endaddress
\email barvinok$\@$umich.edu \endemail
\thanks This research was partially supported by NSF Grants DMS 0400617 and DMS 0856640
and a United States - Israel BSF grant 2006377. \endthanks
 \abstract  Let $R=\left(r_1, \ldots, r_m \right)$ and $C=\left(c_1, \ldots, c_n \right)$ be 
 positive integer vectors such that $r_1 +\ldots + r_m=c_1 +\ldots + c_n$. We consider
 the set $\Sigma(R, C)$ of non-negative $m \times n$ integer matrices (contingency tables)
 with row sums $R$ and column sums $C$ as a finite probability space with the uniform measure. 
 We prove that a random table $D \in \Sigma(R,C)$  is close 
 with high probability to a particular matrix (``typical table'') $Z$ defined as follows. We let 
 $g(x)=(x+1)\ln(x+1)-x \ln x$ for $x \geq 0$ and let $g(X)=\sum_{i,j} g(x_{ij})$ for a non-negative 
 matrix $X=\left(x_{ij}\right)$. Then $g(X)$ is strictly concave and attains its maximum
 on the polytope of non-negative $m \times n$ matrices $X$ with row sums $R$ and column sums
 $C$ at a unique point, which we call the typical table $Z$.
\endabstract
\date November 2009
\enddate
\keywords contingency table, random matrix, transportation polytope
\endkeywords
\subjclass 15A52, 05A16, 60C05, 15A36 \endsubjclass
\endtopmatter

\document

\head 1. Introduction and the main result \endhead

\subhead (1.1) Random contingency tables \endsubhead
Let $R=\left(r_1, \ldots, r_m \right)$ be a positive integer $m$-vector and let 
$C=\left(c_1, \ldots, c_n \right)$ be a positive integer $n$-vector such that 
$$\sum_{i=1}^m r_i =\sum_{j=1}^n c_j =N.$$
A {\it contingency table} with {\it margins} $(R, C)$ is a non-negative integer matrix 
$D=\left(d_{ij}\right)$ with row sums $R$ and column sums $C$:
$$\split \sum_{j=1}^n d_{ij}=r_i \quad &\text{for} \quad i=1, \ldots, m,  \quad  
\sum_{i=1}^m d_{ij} =c_j \quad \text{for} \quad j=1, \ldots, n,\\
&d_{ij} \geq 0 \quad \text{and} \quad d_{ij} \in {\Bbb Z} \quad \text{for all} \quad i,j. \endsplit$$
Let $\Sigma(R, C)$ be the set of all contingency tables with margins $(R, C)$. As is well known,
$\Sigma(R, C)$ is non-empty and finite. Let us consider $\Sigma(R, C)$ as a finite probability 
space endowed with the uniform probability measure. 
In this paper we address the following question:
\smallskip
 Suppose that $D \in \Sigma(R, C)$ 
is chosen at random. What is $D$ likely to look like? 
\smallskip
The problem is interesting in its own right, but the main motivation comes from 
statistics; see \cite{Go63}, \cite{DE85}, \cite{DG95}  and references therein.
A contingency table $D=\left(d_{ij}\right)$ may represent certain statistical data (for example, $d_{ij}$ may be the number of people in a certain sample having the $i$-th hair color and the $j$-th eye color).
One can condition 
on the row and column sums and ask what is special about a particular table 
$D \in \Sigma(R, C)$, considering all tables in $\Sigma(R, C)$ as equiprobable; see \cite{DE85}.
To answer this question we need to know what a random table 
$D \in \Sigma(R,C)$ looks like. Considerable effort 
was invested in finding an efficient (polynomial time) algorithm to {\it sample} a random table 
$D \in \Sigma(R, C)$; see \cite{DG95}, \cite{D+97},
\cite{C+06}. Despite a number of successes, such an algorithm is still at large in many interesting 
situations. In this paper, we do not discuss how to sample a random table but describe instead
what it is likely to look like. 

We prove that a random contingency table $D$ is close in a certain sense 
to some particular non-negative $m \times n$ matrix $Z$, which we call the {\it typical table}.

\subhead (1.2) The typical table \endsubhead Let $\PP(R, C)$ be the set of all $m \times n$ 
non-negative matrices $X=\left(x_{ij}\right)$ with row sums $R$ and column sums $C$:
$$\split \sum_{j=1}^n x_{ij} =r_i \quad &\text{for} \quad i=1, \ldots, m, \quad \sum_{i=1}^m x_{ij} =c_j 
\quad \text{for} \quad j=1, \ldots, n \quad \text{and} \\
&x_{ij} \geq 0 \quad \text{for all} \quad i,j .\endsplit$$
Geometrically, $\PP(R, C)$ is a convex polytope of dimension $(m-1)(n-1)$, known as the 
{\it transportation polytope}. 
Let 
$$g(x)=(x+1) \ln (x+1)-x \ln x \quad \text{for} \quad x \geq 0$$
and let
$$g(X)=\sum_{i,j} g(x_{ij})$$
for a non-negative matrix $X=\left(x_{ij}\right)$.
One can easily check that $g$ is strictly concave and hence achieves a unique maximum 
$Z=\left(z_{ij}\right)$ on $\PP(R, C)$. 
We call $Z$ the {\it typical table} with margins $(R,C)$. 
Since the objective function $g$ is concave, $Z$ can be computed 
efficiently, both in theory and in practice, by existing methods of convex optimization, cf.
\cite{NN94}.

The solution $Z$ to the above optimization problem was first introduced in the author's 
paper \cite{Ba09}. It was given the name of ``typical table'' (perhaps with not enough 
justification) in  \cite{B+08}. 

In this paper, we show that $Z$ indeed captures some typical features of a random table 
$D \in \Sigma(R, C)$.

We prove our main result assuming certain regularity (``smoothness'') of margins.

\subhead (1.3) Smooth margins \endsubhead Let us fix a number $0 < \delta  \leq 1$.
First, we assume that the row sums and column sums are of the same order:
$$\split &{\delta N \over m} \ \leq \ r_i \ \leq \ {N \over \delta m} \quad \text{for} \quad i=1, \ldots, m
\quad \text{and} \\
& {\delta N \over n} \ \leq \ c_j \ \leq \  {N \over \delta n} \quad \text{for} \quad j=1, \ldots, n.
\endsplit \tag1.3.1$$
Second, we assume that the {\it density} of the table is separated from 0:
$${N \over mn}  \geq \delta. \tag1.3.2$$
We say that the margins $(R, C)$ are $\delta$-{\it smooth} if conditions (1.3.1)--(1.3.2) are 
satisfied. This is a modification of the definition from \cite{B+08}. We note that $\delta$-smooth
margins are also $\delta'$-smooth for any $0< \delta' < \delta$. As we remarked (see (1.3.2)),
we are interested in tables with the density separated from 0. For the case of sparse tables, where
$r_i \ll n$ and $c_j \ll m$, see \cite{Ne69}, \cite{GM08} and references therein.

Without loss of generality, we assume that $n \geq m$.
\definition{(1.4) Definitions and notation} Let us choose a non-empty subset of entries of a matrix:
$$S \subset \Bigl\{(i,j): \quad 1 \leq i \leq m, \quad 1 \leq j \leq n \Bigr\}.$$
For an $m \times n$ matrix $A=\left(a_{ij}\right)$ let
$$\sigma_S(A)=\sum_{(i,j) \in S} a_{ij}$$
be the sum of the entries from $S$.

The cardinality of a finite set $X$ is denoted by $|X|$. 
\enddefinition

Now we state our main result.
\proclaim{(1.5) Theorem} Let us fix real numbers $0< \delta \leq 1$ and $\kappa>0$.
Then there exists a positive integer $q = q(\delta, \kappa)$ 
such that the following holds:

Suppose that $(R, C)$ are $\delta$-smooth margins such that $n \geq m \geq q$.

Let 
$$S \subset \left\{(i,j): \quad 1 \leq i \leq m, \quad 1 \leq j \leq n \right\}$$ be a set  such that 
$$|S| \ \geq \ \delta mn,$$
let $Z$ be the typical table with margins $(R, C)$, and let
$$\epsilon =\delta {\ln n \over m^{1/3}}.$$
If $\epsilon \leq 1$ then 
$$\split \Pr \Bigl\{ D \in &\Sigma(R,C): \\ &(1-\epsilon) \sigma_S(Z)  \  \leq \sigma_S(D) \leq  \  (1+\epsilon) \sigma_S(Z)  \Bigr\} \
 \geq \ 1-2n^{-\kappa n}.\endsplit$$
\endproclaim

In  other words, asymptotically, as far as the sum over a positive fraction of entries is concerned,  
a contingency table
$D$ sampled uniformly at random from  the set of contingency tables
with given margins is very likely to be close to the typical table $Z$.

\subhead (1.6) The independence table \endsubhead
In \cite{Go63}, I.J. Good observes that the {\it independence table}  
$$Y=\left(y_{ij}\right), \quad y_{ij}=r_i c_j/N \quad \text{for all} \quad i,j,$$
maximizes the entropy 
$$H(X)=\sum_{i,j} {x_{ij} \over N} \ln {N \over x_{ij}}$$
on the set of all matrices $X=\left(x_{ij}\right)$ in the transportation polytope $\PP(R,C)$.
One may be tempted to think that the independence table $Y$, not the typical table $Z$, 
reflects the structure of a random table $D \in \Sigma(R, C)$.

One can show that $Y=Z$ if and only if all row sums $r_i$ are equal or all column
sums $c_j$ are equal. In fact, particular entries of the matrices $Z$ and $Y$ may demonstrate very 
different behavior even for reasonably looking margins. 
Suppose, for example, that $m=n$, that $r_1=c_1=3n$
and that $r_i=c_i=n$ for $i>1$. Hence $N=3n+n(n-1)=n^2+2n$ and 
for the independence table we have 
$$y_{11}={9n^2 \over n^2 +2n} \leq 9.$$
On the other hand, for the typical table $Z$ the entry $z_{11}$ grows linearly in $n$. 
Indeed, the optimality 
condition for $Z$ (the gradient of $g$ at $Z$ is orthogonal to the affine span of the 
transportation polytope) implies that 
$$\ln \left({z_{ij}+1 \over z_{ij}} \right)=\lambda_i+\mu_j \quad \text{for all} \quad 
i,j$$ 
and some $\lambda_1, \ldots, \lambda_m, \mu_1, \ldots, \mu_n$; see Section 2.3.
By symmetry, we can choose $\lambda_1=\mu_1=\alpha$ and $\lambda_i=\mu_i=\beta$ for $i>1$.
Moreover, we must have $0<\alpha < \beta$.
Since 
$$z_{21}={1 \over e^{\alpha +\beta}-1} > {1 \over e^{2 \beta}-1} =z_{2j} \quad \text{for all} \quad j >1$$
and $r_2=n$, we should have
$$\beta > {\ln 2 \over 2}.$$
Therefore,
$$z_{1j} ={1 \over e^{\alpha + \beta} -1}\ < \ {1 \over e^{\beta}-1} \ < \ {1 \over \sqrt{2}-1} \quad 
\text{for} \quad j >1.$$
Since $r_1=3n$ we must have 
$$z_{11}  \ > \  3n - {n \over \sqrt{2}-1} >0.58 n.$$
Let us show that the independence table $Y$ and the typical table $Z$ may also produce 
different asymptotic behavior of the sums $\sigma_S(Y)$ and $\sigma_S(Z)$ as $m$ and $n$ grow and 
$S$ is a subset of entries consisting of a positive fraction of all entries as in Theorem 1.5.
For that, let us fix some margins $R=\left(r_1, \ldots, r_m \right)$ 
and $C=\left(c_1, \ldots, c_n \right)$ such that $z_{11} \ne y_{11}$. For a positive integer $k$
let us consider the ``cloned'' margins
$$\aligned &R_k=\Bigl(\underbrace{kr_1, \ldots, kr_1}_{\text{ $k$ times}}, \ldots, 
\underbrace{kr_m, \ldots, kr_m}_{\text{$k$ times}} \Bigr) \quad \text{and} \\
&C_k=\Bigl(\underbrace{kc_1, \ldots, kc_1}_{\text{$k$ times}}, \ldots, 
\underbrace{kc_n, \ldots, kr_n}_{\text{$k$ times}} \Bigr). \endaligned \tag1.6.1$$
In particular, tables $D \in \Sigma(R_k, C_k)$ are $km \times kn$ matrices whose 
total sum of entries is equal to $k^2N$, where $N=r_1 + \ldots + r_m =c_1 + \ldots +c_n$.
Let $S=S_k$ be the set of entries in the upper left $k \times k$ corner of a 
matrix from $\Sigma(R_k, C_k)$,
let $Y_k$ be the independence table of margins $(R_k, C_k)$ and let $Z_k$ be the 
typical table of margins $(R_k, C_k)$. It is not hard to show that $\sigma_S(Z_k) = k^2 z_{11}$ and 
$\sigma_S(Y_k)=k^2 y_{11}$, so the ratio between the two sums remains fixed (and not equal to 1)
as $k$ grows. 

It looks plausible that the independence 
table $Y$ is indeed close with high probability to a random table $D \in  \Sigma(R, C)$, if,  instead of the 
uniform distribution in $\Sigma(R,C)$, a table $D=\left(d_{ij}\right)$ is sampled from the 
Fisher-Yates probability measure, where 
$$\Pr(D)=(N!)^{-1} \left(\prod_{i=1}^m r_i! \right) \left(\prod_{j=1}^n c_j! \right) 
\left(\prod_{ij} {1 \over d_{ij}!} \right);$$
see \cite{DG95}.
Compared with the uniform distribution, the Fisher-Yates measure gives less weight to tables with 
large entries.
\bigskip
Let $p, q>0$ be real numbers such that $p+q=1$. Recall that a discrete random variable 
$x$ has {\it geometric distribution} if 
$$\Pr\left\{x=k\right\} =pq^k \quad \text{for} \quad k=0, 1, \ldots$$
We have
$$\EE x ={q \over p}.$$
Consequently, 
$$\text{if} \quad \EE x=z \quad \text{then} \quad p={1 \over 1+z} \quad \text{and} \quad q={z \over 1+z}.$$
The following interpretation of the typical matrix was suggested to the author by J.A. Hartigan; see
\cite{BH09}.
\proclaim{(1.7) Theorem} Let $Z=\left(z_{ij}\right)$ be the $m \times n$ typical table with margins 
$(R, C)$. Let $X=\left(x_{ij}\right)$ be the random $m \times n$ matrix of independent geometric 
random variables $x_{ij}$ such that 
$$\EE x_{ij}=z_{ij} \quad \text{for all} \quad i, j.$$
Then the probability mass function of $X$ is constant on the set $\Sigma(R, C)$ of contingency
tables with margins $(R, C)$, and, moreover,
$$\Pr\left\{X =D \right\} =e^{-g(Z)} \quad \text{for all} \quad D \in \Sigma(R, C),$$
where $g$ is the function defined in Section 1.2.
\endproclaim
In other words, the multivariate geometric distribution $X$ whose expectation is the typical matrix $Z$, when conditioned on the set $\Sigma(R, C)$ of contingency tables, results in the uniform probability 
distribution on $\Sigma(R, C)$. It turns out that for a positive $m \times n$ 
matrix $A$ the value of $g(A)$ is equal to the maximum possible entropy of a random matrix with 
expectation $A$ and 
values in the set 
${\Bbb Z}^{m \times n}_+$ of $m \times n$ non-negative integer matrices.
Such a maximum entropy random matrix is necessarily a matrix 
with independent geometrically distributed entries. Therefore, the distribution of $X$ in Theorem 1.7
can be characterized as the maximum entropy distribution in the class consisting of all probability distributions on 
${\Bbb Z}^{m \times n}_+$  whose expectations 
lie in the affine subspace consisting of the matrices with row sums $R$ and column sums $C$;
see \cite{BH09}.

\subhead (1.8) Possible ramifications and open questions \endsubhead
Theorem 1.7 allows one to interpret Theorem 1.5 as a law 
of large numbers for contingency tables: with respect to sums $\sigma_S(D)$ for sufficiently 
large sets $S$ of entries, a random contingency table $D \in \Sigma(R, C)$ behaves approximately 
as the matrix of independent geometric variables whose expectation is the typical table. 
Similar concentration results can be obtained for other well-behaved functions on contingency 
tables.
One can ask whether the distribution of a {\it particular entry} of a random table 
$D \in \Sigma(R, C)$ is asymptotically geometric, as the dimensions $m$ and $n$ of the 
table grow. For example, does the first entry $d_{11}$ of the table converge in distribution 
to the geometric random variable with expectation $z_{11}$  
when the margins $(R, C)$ are cloned, $(R, C) \longmapsto \left(R_k, C_k \right)$, as in (1.6.1)? 

Let us fix a subset 
$$W \subset \Bigl\{(i,j): \quad i=1, \ldots, m; \ j=1, \ldots, n \Bigr\}.$$
Let us consider the set $\Sigma(R, C; W)$ of $m \times n$ non-negative integer matrices 
$D=\left(d_{ij}\right)$ with row sums $R$, column sums $C$ and such that $d_{ij}=0$ for 
$(i,j) \notin W$. Assuming that $\Sigma(R, C; W)$ is non-empty, we can consider $\Sigma(R, C; W)$
as a finite probability space with the uniform measure and ask what a random table 
$D \in \Sigma(R,C; W)$ looks like. 

As above, we define the typical table $Z$ as the unique 
maximum of $g(X)$ on the polytope of non-negative matrices $X=\left(x_{ij} \right)$ with 
row sums $R$, column sums $C$ and such that $x_{ij}=0$ for $(i,j) \notin W$. One can 
prove versions of Theorem 1.5 and Theorem 1.7 in this more general context for subsets $S \subset W$.
However, it appears that for Theorem 1.5 one has to assume, additionally, that there are no  
too large or too small values among the entries $z_{ij}$ of the typical table
$Z=\left(z_{ij}\right)$, cf. the example in Section 1.6. In our case, when $W$ is the set 
of all pairs $(i,j)$, Lemma 2.4 ensures that the entries $z_{ij}$ are not too small while 
Lemma 3.3 ensures that they are not too large.

In \cite{Ba08} another variation of the problem is considered: what if we require 
$d_{ij} \in \{0,1\}$ for all $i,j$. It turns out that a random $D$ is close to a particular matrix maximizing 
the sum of entropies of the entries among all matrices with row sums $R$, column sums $C$ and entries between 
0 and 1.

\bigskip
In the rest of the paper, we prove Theorem 1.5.

In Section 2, we recall the main results of \cite{Ba09} connecting the typical table $Z$ with an 
asymptotic estimate for the number $|\Sigma(R,C)|$ of tables and also prove Theorem 1.7.

In Section 3, we prove Theorem 1.5 under the additional assumption that the total sum $N$ of 
the entries is bounded 
by a polynomial in $m$ and $n$.

In Section 4, we complete the proof of Theorem 1.5.

\head 2. Preliminaries: an asymptotic formula for the number of tables \endhead

In \cite{Ba09}, the following result was proved; see Theorem 1.1 there.

\proclaim{(2.1) Theorem} Let $R=\left(r_1, \ldots, r_m \right)$ and $C=\left(c_1, \ldots, c_n \right)$
be positive integer vectors such that $r_1 + \ldots + r_m =c_1 + \ldots + c_n=N$.
Let us define a function
$$\split F(&\xx, \yy)=\left(\prod_{i=1}^m x_i^{-r_i} \right) \left( \prod_{j=1}^n y_j^{-c_j} \right)
\left( \prod_{i, j} {1 \over 1 - x_i y_j} \right) \\ &\quad \text{for} \quad
\xx=\left(x_1, \ldots, x_m \right) \quad \text{and} \quad \yy=\left(y_1, \ldots, y_n \right).
\endsplit$$
Then $F(\xx, \yy)$ attains its infimum
$$\rho(R,C)=\min  \Sb 0 < x_1, \ldots, x_m <1 \\ 0< y_1, \ldots, y_n <1 \endSb F(\xx, \yy)$$
on the open cube $0<x_i, y_j<1$ and 
for the number $|\Sigma(R,C)|$ of non-negative integer $m \times n$ matrices with row 
sums $R$ and column sums $C$ we have
$$\rho(R, C) \ \geq \ |\Sigma(R,C)| \ \geq \ N^{-\gamma(m+n)} \rho(R,C),$$
where $\gamma >0$ is an absolute constant.
\endproclaim
{\hfill \hfill \hfill} \qed

As is remarked in \cite{Ba09}, the substitution $x_i=e^{-s_i}$, $y_j=e^{-t_j}$ transforms 
$\ln F(\xx, \yy)$ into a convex function 
$$\split G(&\sss, \ttt)=\sum_{i=1}^m r_i s_i + \sum_{j=1}^n c_j t_j -\sum_{i,j} \ln \left(1-e^{-s_i -t_j}\right)
\\ &\quad \text{for} \quad \sss=\left(s_1, \ldots, s_m\right) \quad \text{and} \quad
\ttt=\left(t_1, \ldots, t_n \right) \endsplit$$
on the positive orthant ${\Bbb R}^m_+ \times {\Bbb R}^n_+$. 
It turns out that the typical table $Z$ is the solution to the problem that is convex dual to 
the problem of minimizing $G$. The following result was proved in \cite{Ba09}; see Lemma 1.4
there.

\proclaim{(2.2) Lemma} Let $\PP=\PP(R,C)$ be the polytope of $m \times n$ non-negative matrices
$X=\left(x_{ij}\right)$ with row sums $R$ and column sums $C$ and let $Z \in \PP(R,C)$ 
be the typical table; see Section 1.2.

Then one can write $Z=\left(z_{ij}\right)$,
$$z_{ij}={\xi_i \eta_j \over 1-\xi_i \eta_j} \quad \text{for all} \quad i,j$$
and some $0< \xi_1, \ldots, \xi_m; \eta_1, \ldots, \eta_n <1$ such that the minimum $\rho(R,C)$ 
of the function  $F(\xx, \yy)$ in Theorem 2.1 is attained at 
$x^{\ast}=\left(\xi_1, \ldots, \xi_m \right)$ and \break $\yy^{\ast}=\left(\eta_1, \ldots, \eta_n \right)$:
$$F\left(\xx^{\ast}, \yy^{\ast}\right) =\rho(R,C)=\min \Sb 0 < x_1, \ldots, x_m <1 \\ 0< y_1, \ldots, y_n < 1 \endSb
F(\xx, \yy).$$
Moreover,
$$\rho(R,C)=\exp\left\{g(Z)\right\}.$$
\endproclaim
{\hfill \hfill \hfill} \qed

Theorem 1.7 is a particular case of a more general result proved in \cite{BH09}. Nevertheless, 
we present the proof of Theorem 1.7 here for completeness and 
since some elements of the proof will be recycled later.

\subhead (2.3) Proof of Theorem 1.7 \endsubhead 
From Lemma 2.2, we have $z_{ij} >0$ for all $i,j$. Since $Z$ lies in the relative interior of the transportation polytope
$\PP(R, C)$, the gradient of $g$ at $Z$ must be orthogonal to the subspace of $m \times n$ 
matrices with row and column sums equal to 0.
Therefore,
$$\ln \left({z_{ij}+ 1\over z_{ij}} \right) =\lambda_i +\mu_j  \quad \text{for all} \quad i, j \tag2.3.1$$
and some $\lambda_1, \ldots, \lambda_m$ and $\mu_1, \ldots, \mu_n$. 

For the geometric random variables $x_{ij}$ we have 
$$\Pr\Bigl\{ x_{ij}=d_{ij} \Bigr\}=p_{ij}q_{ij}^{d_{ij}} =\left({1 \over 1+z_{ij}}\right)
 \left({z_{ij} \over 1+z_{ij}} \right)^{d_{ij}}$$
Using (2.3.1), for $D \in \Sigma(R, C)$, $D=\left(d_{ij}\right)$, we obtain
$$\split \Pr\bigl\{X=D\bigr\} 
=&\left( \prod_{i,j} {1 \over 1+z_{ij}} \right)  \prod_{i, j} \left({z_{ij} \over 1+z_{ij}}\right)^{d_{ij}} \\
=&\left( \prod_{i,j} {1 \over 1+z_{ij}} \right)  \prod_{i, j} e^{-(\lambda_i +\mu_j)d_{ij}} \\
=&\left( \prod_{i,j} {1 \over 1+z_{ij}} \right)  \left(\prod_{i=1}^m e^{-\lambda_i r_i} \right) 
\left( \prod_{j=1}^n e^{-\mu_j c_j} \right).
\endsplit $$
Also, 
$$\split e^{-g(Z)} =&\prod_{i,j} {z_{ij}^{z_{ij}} \over \left(1+z_{ij}\right)^{z_{ij} +1}}\\
 =&\left(\prod_{i,j} {1 \over 1+z_{ij}}\right) \prod_{i,j} \left({ z_{ij} \over 1+ z_{ij}} \right)^{z_{ij}} \\
 =&\left(\prod_{i,j} {1 \over 1+z_{ij}}\right) \prod_{i,j} e^{-(\lambda_i +\mu_j)z_{ij}}\\
 =&\left( \prod_{i,j} {1 \over 1+z_{ij}} \right)  \left(\prod_{i=1}^m e^{-\lambda_i r_i} \right) 
\left( \prod_{j=1}^n e^{-\mu_j c_j} \right),
 \endsplit$$
which completes the proof. 
{\hfill \hfill \hfill} \qed

We will need a lower bound for the entries of the typical table $Z=\left(z_{ij}\right)$ proved in
\cite{B+08}; see Theorem 3.3 there.
\proclaim{(2.4) Lemma} Let
$$\split &r_+=\max \Sb i=1, \ldots, m \endSb r_i, \quad r_-=\min \Sb i=1, \ldots, m \endSb r_i 
\quad \text{and} \\
 &c_+=\max \Sb j=1, \ldots, n \endSb c_j, \quad c_-=\min \Sb j=1, \ldots, n \endSb c_j. \endsplit$$
Let $Z=\left(z_{ij}\right)$ be the typical table with margins $(R, C)$. Then
$$z_{ij} \ \geq \ {r_- c_- \over r_+ m} \quad \text{and} \quad z_{ij} \ \geq \ {c_- r_- \over c_+ n} 
\quad \text{for all} \quad i, j.$$
\endproclaim
{\hfill \hfill \hfill} \qed

\proclaim{(2.5) Corollary} Let $Z=\left(z_{ij}\right)$ be the typical table of $\delta$-smooth 
margins $(R, C)$. Then 
$$z_{ij} \ \geq \ {\delta^3 N \over mn} \quad \text{for all} \quad i, j.$$
\endproclaim
\demo{Proof} In Lemma 2.4, we have 
$$r_- \ \geq \ {\delta N \over m},  \quad c_- \ \geq \ {\delta N \over n} \quad \text{and} \quad 
r_+ \ \leq \ {N \over \delta m},$$
and the result follows.
{\hfill \hfill \hfill} \qed
\enddemo

\head 3. Proof of Theorem 1.5 assuming that $N$ is polynomially bounded  \endhead

In this section we prove Theorem 1.5 under the additional assumption that the total sum $N$
of entries is bounded
by a polynomial in $m$ and $n$, specifically that $N \leq (mn)^{1/\delta}$. We use Theorem 1.7.
We start with a standard large deviation inequality.

\proclaim{(3.1) Lemma} Let $X=\left(x_{ij} \right)$ be the $m \times n$ matrix of independent 
geometric random variables $x_{ij}$ such that $\EE X=Z$, $Z=\left(z_{ij}\right)$. 
Let 
$$S \subset \Bigl\{ (i, j): \quad 1 \leq i \leq m, \quad 1 \leq j \leq n \Bigr\}$$
be a non-empty set.
Recall that
$$\sigma_S(X)=\sum_{(i, j) \in S} x_{ij}, \quad \sigma_S(Z)= \sum_{(i, j) \in S} z_{ij}$$
and let us denote
$$\nu_S(Z) =\sum_{(i, j) \in S} z_{ij}^2.$$
Then 
\roster
\item For any real $a$ and for any $0< t \leq 2$, we have 
$$\Pr\bigl\{\sigma_S(X) \ \leq \ - a + \sigma_S(Z) \bigr\} \ \leq \ \exp\left\{ - t a + {t^2 \over 2}
 \bigl(\sigma_S(Z) + \nu_S(Z) \bigr) \right\}.$$
\item For any real $a$ and for any $0 < t \leq \min\bigl\{1/3,\ 1/2z_{ij}: \ (i,j) \in S\bigr\}$,
we have 
$$\Pr\bigl\{\sigma_S(X) \ \geq \ a + \sigma_S(Z) \bigr\} \ \leq \ \exp\Bigl\{ - t a + 2t^2
 \bigl(\sigma_S(Z) + \nu_S(Z) \bigr) \Bigr\}.$$
\endroster
\endproclaim
\demo{Proof}
We use the Laplace transform method; see, for example, Section 1.6 of \cite{Le01}. 
To prove Part (1), for any $t>0$ we compute
$$\EE e^{-t\sigma_S(X)} =\prod_{(i, j) \in S} \EE e^{-t x_{ij}}= \prod_{(i, j) \in S} {p_{ij} \over 1-e^{-t} q_{ij}},$$
where
$$\Pr\bigl\{ x_{ij} =k \bigr\} = p_{ij} q_{ij}^k \quad \text{for} \quad k=0, 1, \ldots$$
Using the fact that $e^{-t} \leq 1-t+t^2/2$ for $t \geq 0$, we obtain 
$$\EE e^{-t \sigma_S(X)} \  \leq \ \prod_{(i, j) \in S} {p_{ij} \over p_{ij} + (t-t^2/2) q_{ij}}= \prod_{(i, j) \in S} 
{1 \over 1 + (t-t^2/2) z_{ij}}.$$ 
Using the fact that $t-t^2/2 \geq 0$ for $0\leq t \leq 2$ and that $\ln (1+x) \geq x -x^2/2$ for $x \geq 0$, we obtain 
$$\split \EE e^{-t\sigma_S(X)} \ \leq \ &\exp\left\{ -\sum_{(i, j) \in S} \ln \bigl(1 +\left(t-t^2/2\right) z_{ij} 
\bigr) \right\} \\
 \leq \ &\exp\left\{ - \sum_{(i, j) \in S } (t-t^2/2) z_{ij} + {1 \over 2} \sum_{(i, j) \in S}  (t-t^2/2)^2 z_{ij}^2 \right\}
 \\ \leq \ & \exp\left\{ - t \sigma_S(Z) + {t^2 \over 2} \bigl(\sigma_S(Z) + \nu_S(Z) \bigr) \right\}.  \endsplit$$
 Then 
 $$\split \Pr \bigl\{\sigma_S(X) \  \leq \ - a + \sigma_S(Z)  \bigr\} \ =\ &\Pr \bigl\{ -t \sigma_S(X) \  \geq  \ ta - t \sigma_S(Z) \bigr\} \\ = \ 
& \Pr \left\{ e^{-t\sigma_S(X)} \  \geq  \ e^{ ta - t \sigma_S(Z)} \right\} \\
\  \leq \ &e^{-ta + t \sigma_S(Z)} \EE e^{-t\sigma_S(X)} \\
\  \leq \ &\exp\left\{-ta + {t^2 \over 2} \bigl(\sigma_S(Z) + \nu_S(Z) \bigr)  \right\}. \endsplit$$

To prove Part (2), we observe that $e^t < 1+t+t^2$ for all $0 < t \leq 1$. 
Therefore, for $0< t \leq \min\bigl\{1/3, \ 1/2z_{ij}: \ (i,j) \bigr\}$, we have 
$$e^t \ < \ 1+2t \ \leq \ {1+z_{ij} \over z_{ij}} ={1 \over q_{ij}}$$ 
and hence
$$\split \EE e^{t \sigma_S(X)} =&\prod_{(i, j) \in S} \EE e^{t x_{ij}} =
\prod_{(i, j) \in S} {p_{ij} \over 1-e^t q_{ij}} \\
\leq &\prod_{(i, j) \in S} {p_{ij} \over p_{ij} -(t+t^2) q_{ij}} = \prod_{(i, j) \in S}
{1 \over 1-(t+t^2) z_{ij}}.
\endsplit $$
Since $t \leq 1/3$ we have $t+t^2 \leq (4/3) t$ and hence $(t+t^2) z_{ij} \leq 2/3$.
Using the fact that $\ln (1-x) \geq -x -x^2$ for $0 \leq x \leq 2/3$, we obtain
$$\split \EE e^{t \sigma_S(X)} \ \leq \ &\exp\left\{ - \sum_{(i, j) \in S} \ln \left(1 -(t+t^2) z_{ij} \right)
\right\} \\ 
\leq \ &\exp\left\{ \sum_{(i, j) \in S} \left(t + t^2 \right) z_{ij} + \sum_{(i, j) \in S} \left( t+ t^2 \right)^2 
z_{ij}^2 \right\}  \\
\leq \ &\exp\Bigl\{ t \sigma_S(Z)  +2 t^2 \left(\sigma_S(Z) + \nu_S(Z) \right) \Bigr\}.\endsplit$$
Therefore,
 $$\split \Pr \bigl\{\sigma_S(X) \  \geq \  a + \sigma_S(Z)  \bigr\} \ =\ &\Pr \bigl\{ t \sigma_S(X) \  \geq  \ ta + t \sigma_S(Z) \bigr\} \\ = \ 
& \Pr \left\{ e^{t\sigma_S(X)} \  \geq  \ e^{ ta + t \sigma_S(Z)} \right\} \\
\  \leq \ &e^{-ta - t \sigma_S(Z)} \EE e^{t\sigma_S(X)} \\
\  \leq \ &\exp\Bigl\{-ta + 2t^2\bigl(\sigma_S(Z) + \nu_S(Z) \bigr)  \Bigr\}. \endsplit$$
 {\hfill \hfill \hfill} \qed
\enddemo
One can observe that $\sigma_S(Z) + \nu_S(Z)$ is the variance of $\sigma_S(X)$. 
\proclaim{(3.2) Corollary} Let $(R, C)$ be $\delta$-smooth margins with the typical 
table $Z=\left(z_{ij}\right)$ and let $X=\left(x_{ij}\right)$ be the matrix of independent geometric 
variables such that $\EE X=Z$. Suppose that 
$$z_{ij} \ \leq \ {\alpha N \over mn} \quad \text{for all} \quad (i, j) \in S$$
and some $\alpha \geq 1$.
Then
\roster
\item For any $0 < \epsilon < 1$ we have 
$$\Pr \bigl\{ \sigma_S(X)  \ \leq \ (1-\epsilon) \sigma_S(Z) \bigr\} \ \leq \ \exp\left\{ -{\epsilon^2 \delta^4 |S| \over 2 + 2 \delta \alpha} \right\}.$$
\item For any $0 < \epsilon < 1$ we have 
$$\Pr\bigl\{\sigma_S(X) \ \geq \ (1+\epsilon) \sigma_S(Z) \bigr\} 
\ \leq \ \exp\left\{ -{\epsilon^2 \delta^4 |S| \over 8 + 8 \delta \alpha} \right\}.$$
\endroster
\endproclaim
\demo{Proof} 
Choosing
$$a=\epsilon \sigma_S(Z) \quad \text{and} \quad t=
{\epsilon \sigma_S(Z) \over \sigma_S(Z) + \nu_S(Z)}$$
in Part (1) of  Lemma 3.1, we obtain 
$$\Pr \bigl\{ \sigma_S(X)  \ \leq \ (1-\epsilon) \sigma_S(Z) \bigr\} \ \leq \ 
\exp\left\{-{\epsilon^2 \sigma_S^2(Z) \over 2(\sigma_S(Z) + \nu_S(Z))} 
\right\}. \tag3.2.1$$
Furthermore,
$$\nu_S(Z)=\sum_{(i, j) \in S} z_{ij}^2 \ \leq \ {\alpha N \over mn}  \sum_{(i, j) \in S} z_{ij} \ = \ 
{\alpha N \over mn}  \sigma_S(Z). \tag3.2.2$$
By Corollary 2.5,
$$\delta_S(Z) \ \geq \ |S| {\delta^3 N \over mn}. \tag3.2.3$$
We recall that 
$${N \over mn} \ \geq \ \delta.\tag3.2.4$$
Summarizing (3.2.1)--(3.2.4), we get 
$$\split \Pr \bigl\{ \sigma_S(X)  \ \leq \ (1-\epsilon) \sigma_S(Z) \bigr\} \ \leq \ 
&\exp\left\{ -{\epsilon^2 \sigma_S(Z) mn \over 2\left(mn + \alpha N \right)} \right\} \\
\leq \  & \exp\left\{ -{\epsilon^2 |S| \delta^3 N \over 2\left(mn + \alpha N\right)} \right\} \\
\leq \ & \exp\left\{ -{\epsilon^2 |S| \delta^3  \over 2 \left(mn/N + \alpha \right)} \right\} \\
\leq \ & \exp\left\{ -{\epsilon^2 \delta^4 |S| \over 2 + 2 \delta \alpha} \right\} \endsplit$$
and Part (1) follows. 

Let us choose $a=\epsilon \sigma_S(Z)$ in Part (2) of Lemma 3.1.
Let 
$$t_0 ={\epsilon \sigma_S(Z) \over 4 \left( \sigma_S(Z) + \nu_S(Z)\right)}.$$
Clearly, $t_0 \leq 1/4 < 1/3$. If $t_0 < mn/2\alpha N$, we choose $t=t_0$ and 
if $t_0  \geq mn/2 \alpha N$, we choose $t=mn/2 \alpha N$ in Part (2) of Lemma 3.1.
Hence if $t_0 < mn/2 \alpha N$, we obtain as above in Part (1)
$$\aligned \Pr\bigl\{ \sigma_S(X) \ \geq \ (1+\epsilon) \sigma_S(Z) \bigr\} \ \leq \ 
&\exp\left\{-{\epsilon^2 \sigma_S^2(Z) \over 8 \left(\sigma_S(Z) + \nu_S(Z)\right)} \right\} \\
\leq \ &\exp\left\{ -{\epsilon^2 \delta^4 |S| \over 8 + 8 \delta \alpha} \right\}. \endaligned \tag3.2.5$$
If $t_0 \geq mn/2 \alpha N$ then 
$$\sigma_S(Z) +\nu_S(Z) \ \leq \ {\epsilon \sigma_S(Z) \alpha N \over 2 mn}.$$  
Therefore, choosing $t=mn/2\alpha N$ in Part (2) of Lemma 3.1, we obtain 
$$\split \Pr\bigl\{ \sigma_S(X) \ \geq \ (1+\epsilon) \sigma_S(Z) \bigr\} \ \leq \ 
&\exp\left\{-{\epsilon mn \sigma_S(Z) \over 2 \alpha N} +{m^2n^2 \left( \sigma_S(Z) +\nu_S(Z) \right) \over 2 \alpha^2 N^2}   \right\} \\
\leq \ &\exp\left\{ - {\epsilon \sigma_S(Z) mn \over 4 \alpha N} \right\}.  \endsplit  $$
Using (3.2.3), we obtain 
$$\Pr \bigl\{ \sigma_S(X) \ \geq \ (1+\epsilon) \sigma_S (Z) \bigr\} \ \leq \ 
\exp\left\{ -{ \epsilon \delta^3 |S| \over 4 \alpha} \right\}. \tag3.2.6$$
Comparing (3.2.5) and (3.2.6), we complete the proof.
{\hfill \hfill \hfill} \qed
\enddemo

Now we can prove the following weaker version of Theorem 1.5.
\proclaim{(3.3) Proposition} Let us fix real numbers $0 < \delta \leq 1$ and $\kappa >0$. Then 
there exists a positive integer $q=q(\delta, \kappa)$ such that the following holds:

Suppose that $(R, C)$ are $\delta$-smooth margins such that $n \geq m \geq q$ and let $Z=\left(z_{ij}\right)$ be the typical table with margins $(R, C)$.
Let 
$$S \subset \bigl\{(i, j): \quad 1 \leq i \leq m, \quad 1 \leq j \leq n \bigr\}$$
be a set such that 
$$|S| \geq \delta mn$$ and suppose that the entries $z_{ij}$ of the typical table satisfy the 
inequalities
$$z_{ij} \ \leq \ {\alpha N \over mn} \quad \text{for} \quad \alpha=2 \delta^{-1} m^{1/3}$$ 
and all $(i, j ) \in S$. 

Suppose further that for the  total sum $N$ of entries we have
$$N \ \leq \ (mn)^{1/\delta}.$$
Let 
$$\epsilon={\delta \ln n \over m^{1/3}}.$$ 
If $\epsilon \leq 1$, we have 
$$ \split &\Pr\bigl\{D \in \Sigma(R, C): \ \sigma_S(D) \ \leq \ (1-\epsilon) \sigma_S(Z) \bigr\}  \ \leq \ 
n^{-\kappa n} \quad \text{and} \\
&\Pr\bigl\{D \in \Sigma(R, C): \ \sigma_S(D) \ \geq \ (1+\epsilon) \sigma_S(Z) \bigr\} \ \leq \
n^{-\kappa n}. \endsplit$$
\endproclaim
\demo{Proof} Let $X=\left(x_{ij}\right)$ be the $m \times n$ matrix of independent geometric random variables 
$x_{ij}$ such that $\EE X= Z$. By Theorem 1.7, the distribution of $X$ conditioned on 
$X \in \Sigma(R, C)$ is uniform and hence
$$\split &\Pr \bigl\{D \in \Sigma(R, C): \  \sigma_S(D) \leq (1-\epsilon) \sigma_S(Z) \bigr\} \\
&\qquad \qquad= {\Pr \bigl\{X: \  \sigma_S(X) \leq (1-\epsilon) \sigma_S(Z) \quad
\text{and} \quad X \in \Sigma(R, C)  \bigr\} \over \Pr\bigl\{X: \ X \in \Sigma(R, C) \bigr\}}.
\endsplit$$
Similarly, 
$$\split &\Pr \bigl\{D \in \Sigma(R, C): \  \sigma_S(D) \geq (1+\epsilon) \sigma_S(Z) \bigr\} \\
&\qquad \qquad= {\Pr \bigl\{X: \  \sigma_S(X) \geq (1+\epsilon) \sigma_S(Z) \quad
\text{and} \quad X \in \Sigma(R, C)  \bigr\} \over \Pr\bigl\{X: \ X \in \Sigma(R, C) \bigr\}}.
\endsplit$$

By Theorem 1.7, Lemma 2.2 and Theorem 2.1 we get  
$$\Pr \bigl\{ X \in \Sigma(R, C) \bigr\} \ = \ e^{-g(Z)} \left| \Sigma(R, C) \right| \ \geq \ N^{-\gamma(m+n)}$$
for some absolute constant $\gamma >0$. Since $N \leq (mn)^{1/\delta}$, we obtain
$$\aligned &\Pr \bigl\{D \in \Sigma(R, C): \ \sigma_S(D) \leq (1-\epsilon) \sigma_S(Z) \bigr\} \\
&\qquad \qquad \leq \ (mn)^{\gamma_1(m+n)}
  \Pr \bigl\{X:  \ \sigma_S(X) \leq (1-\epsilon) \sigma_S (Z) \bigr\} \endaligned$$
and similarly
$$\aligned &\Pr \bigl\{D \in \Sigma(R, C): \ \sigma_S(D) \geq (1+\epsilon) \sigma_S(Z) \bigr\} \\
&\qquad \qquad \leq \ (mn)^{\gamma_1(m+n)}
  \Pr \bigl\{X:  \ \sigma_S(X) \geq (1+\epsilon) \sigma_S (Z) \bigr\} \endaligned$$
for some constant $\gamma_1=\gamma(\delta)>0$.
By Part  (1) of Corollary 3.2, 
$$ \Pr \bigl\{X:\  \sigma_S(X) \leq (1-\epsilon) \sigma_S (Z) \bigr\} \ \leq \ 
\exp\left\{-{\delta^7 mn  \ln^2 n   \over m^{2/3} (2 + 4 m^{1/3})}\right\},$$
while by Part (2) of Corollary 3.2
$$ \Pr \bigl\{X:\  \sigma_S(X) \geq (1+\epsilon) \sigma_S (Z) \bigr\} \ \leq \ 
\exp\left\{-{\delta^7 mn  \ln^2 n   \over m^{2/3} (8 + 16 m^{1/3})}\right\},$$
and the result follows.
{\hfill \hfill \hfill} \qed 
\enddemo 

Next, we prove that large entries of the typical table $Z$ belong to a small number of rows.

\proclaim{(3.4) Lemma} Let $(R, C)$ be $\delta$-smooth margins and let 
$Z=\left(z_{ij}\right)$ be the $m \times n$ typical table with margins $(R, C)$.
 Let $\alpha \ \geq \ 2mn/N$ be a real number.
Let
$$I=\Bigl\{i: \quad z_{ij} \ \geq \  {\alpha N \over mn} \quad \text{for some} \quad j \Bigr\}.$$
Then
$$|I| \ \leq \ {4m \over \delta \alpha}.$$
\endproclaim
\demo{Proof} By (2.3.1), we can write
$$\ln \left({z_{ij}+1 \over z_{ij}} \right) =\lambda_i +\mu_j \quad \text{for all} \quad i, j$$
and some $\lambda_1, \ldots, \lambda_m$ and $\mu_1, \ldots, \mu_n$. Since 
$\lambda_i +\mu_j >0$ for all $i$ and $j$, without loss of generality we may assume that 
$\lambda_1, \ldots, \lambda_m$ and $\mu_1, \ldots, \mu_n$ are positive.

Let
$$I_0=\Bigl\{i: \quad \lambda_i\ \leq \ {mn \over \alpha N} \Bigr\} \quad \text{and} \quad 
J_0 =\Bigl\{ j: \quad \mu_j \ \leq \ {mn \over \alpha N} \Bigr\}.$$
If $i \in I$ then for some $j$ we have 
$${mn \over \alpha N} \ \geq \ 
 {1 \over z_{ij}}\ \geq \ \ln \left({z_{ij} +1 \over z_{ij}} \right) \ \geq \ \lambda_i 
$$
and therefore $I \subset I_0$.
Similarly, if $z_{ij} \geq \alpha N/mn$ for some $i$ then $j \in J_0$. 
Hence without loss of generality, we may assume that $J_0 \ne \emptyset$.

Let us fix a $j_0 \in J_0$. Then for any $i \in I_0$ we have 
$$\ln \left({z_{ij_0} +1 \over z_{ij_0}} \right) \ \leq \ {2 mn \over \alpha N}.$$
Hence for all $i \in I_0$ we have 
$${1 \over z_{ij_0}} \ \leq \ \exp\left\{{2mn \over \alpha N} \right\} -1 \ \leq \ {4 mn \over \alpha N}$$
(using the fact that $e^x \leq 1+2x$ for $0 \leq x \leq 1$).
Hence
$$z_{ij_0} \ \geq \ {\alpha N \over 4mn} \quad \text{for} \quad i \in I_0.$$
Since
$$\sum_{i=1}^m z_{ij_0} =c_{j_0} \ \leq \ {N \over \delta n},$$
we conclude that 
$$|I|\ \leq \ |I_0| \ \leq \ {4 c_{j_0}  mn \over \alpha N} \ \leq \ {4 m \over \delta \alpha}.$$
{\hfill \hfill \hfill} \qed 
\enddemo

Finally, we prove the main result of this section.
\proclaim{(3.5) Proposition} In Theorem 1.5 assume, additionally, that $N \leq (mn)^{1/\delta}$
(equivalently, drop the upper bound assumption for $z_{ij}$ in Proposition 3.3).
Then the conclusion of Theorem 1.5 holds (equivalently, the conclusion of Proposition 3.3 holds).
\endproclaim 

\demo{Proof}
Let us choose 
$$\alpha=2 \delta^{-1} m^{1/3}$$
and let 
$$I=\left\{i: \quad z_{ij} \geq {\alpha N \over mn} \quad \text{for some} \quad j \right\}.$$
Since $N/mn \geq  \delta$, we have $\alpha \geq 2mn/N$ and by Lemma 3.4
we have 
$$|I|  \ \leq \ 2 m^{2/3}.$$
Let 
$$S_0=\bigl\{(i, j)\in S: \quad i \notin I \bigr\}.$$
Then 
$$\left|S \setminus S_0\right| \ \leq \ n |I| \ \leq \ 2 n m^{2/3},$$ 
and hence for $\delta_0=\delta/2$ and $m$ sufficiently large, $n \geq m  \geq q(\delta)$, we have 
$$|S_0| \ \geq \ \delta_0 mn.$$
Furthermore, we have 
$$\sigma_{S \setminus S_0}(D), \ \sigma_{S \setminus S_0} (Z) \ \leq \ \sum_{i \in I} r_i \ \leq  \ |I| {N \over \delta m} 
\ \leq \ {2N \over \delta m^{1/3}}.$$
On the other hand, by Corollary 2.5, we have 
$$\sigma_S (Z) \ \geq \ |S| {\delta^3 N \over mn} \ \geq \ \delta^4 N.$$
Therefore,
$$\sigma_{S_0}(Z)\ = \ \sigma_S(Z)-\sigma_{S\setminus S_0}(Z) \ \geq \ \left(1-{2\over \delta^5 m^{1/3}}\right) \sigma_S(Z) \tag3.5.1$$
and, similarly,
$$\sigma_S(Z) -\sigma_{S \setminus S_0}(D) \ \geq \ \left(1-{2\over \delta^5 m^{1/3}}\right) \sigma_S(Z). \tag3.5.2$$
We have
$$\split &\Pr \left\{D \in \Sigma(R, C):\  \sigma_S(D) \leq (1-\epsilon) \sigma_S(Z) \right\} \\
  &\qquad \qquad \leq \
\Pr\left\{D \in \Sigma(R, C): \ \sigma_{S_0}(D) \leq (1-\epsilon) \sigma_S(Z) \right\}. \endsplit$$
By (3.5.1) we obtain
 $$\split (1-\epsilon) \sigma_S(Z) \ = \ &\left(1 - {\delta \ln n \over m^{1/3}}\right) \sigma_S(Z) 
 \ \leq \ \left(1-\epsilon_0 \right) \left(1-{2 \over \delta^5 m^{1/3}} \right) \sigma_S(Z) \\
 \ \leq \ &\left(1-\epsilon_0\right) \sigma_{S_0}(Z), \quad \text{where}\\
 &\qquad \ \epsilon_0 = {\delta \ln n \over 2 m^{1/3}}, \endsplit$$
 and $m$ is sufficiently large,  $n \geq m \geq q(\delta)$.
 
 Applying Proposition 3.3  with $S_0 \subset S$ and $\delta_0=\delta/2$, we conclude that if $m$ is sufficiently large,
$n \geq m \geq  q(\delta, \kappa)$, we have 
$$\split 
&\Pr\bigl\{D \in \Sigma(R, C): \ \sigma_{S_0}(D) \leq (1-\epsilon) \sigma_S(Z) \bigr\} \\
&\qquad \qquad \leq \
\Pr\bigl\{D \in \Sigma(R, C): \ \sigma_{S_0}(D) \leq \left(1-\epsilon_0\right) \sigma_{S_0}(Z) \bigr\}
\ \leq \ n^{-\kappa n}. \endsplit$$ 
Similarly, we have 
$$\split &\Pr \bigl\{D \in \Sigma(R, C):\  \sigma_S(D) \geq (1+\epsilon) \sigma_S(Z) \bigr\} \\
  &\qquad \qquad = \
\Pr\bigl\{D \in \Sigma(R, C): \ \sigma_{S_0}(D) \geq (1+\epsilon) \sigma_S(Z) -
\sigma_{S \setminus S_0} (D)  \bigr\} \\ 
 &\qquad \qquad \leq \
\Pr\bigl\{D \in \Sigma(R, C): \ \sigma_{S_0}(D) \geq (1+\epsilon) \bigl( \sigma_S(Z) -
\sigma_{S \setminus S_0} (D) \bigr)  \bigr\}. \endsplit$$
By (3.5.2) we obtain 
$$\split (1+\epsilon) \bigl( \sigma_S(Z) - \sigma_{S \setminus S_0}(D) \bigr) \ \geq \ 
&(1+ \epsilon) \left(1 - {2 \over \delta^5 m^{1/3}} \right) \sigma_S(Z) \\ \geq \ 
&\left(1+ \epsilon_0\right) \sigma_{S_0}(Z), \quad \text{where} \\
&\epsilon_0={\delta \ln n \over 2 m^{1/3}}, \endsplit$$
and $m$ is sufficiently large, $n \geq m \geq q(\delta)$.

 Applying Proposition 3.3  with $S_0 \subset S$ and $\delta_0=\delta/2$, we conclude that if $m$ is sufficiently large,
$n \geq m \geq  q(\delta, \kappa)$, we have 
$$\split 
&\Pr\bigl\{D \in \Sigma(R, C): \ \sigma_{S_0}(D) \geq \left(1+\epsilon\right) 
\left(\sigma_S(Z) - \sigma_{S \setminus S_0}(D) \right) \bigr\} \\
&\qquad \qquad \leq \
\Pr\bigl\{D \in \Sigma(R, C): \ \sigma_{S_0}(D) \geq \left(1+\epsilon_0\right) \sigma_{S_0}(Z) \bigr\}
\ \leq \ n^{-\kappa }\endsplit$$
and the result follows. 
{\hfill \hfill \hfill} \qed
\enddemo

\head 4.  Proof of Theorem 1.5 \endhead

It remains to prove Theorem 1.5 in the case of a large (superpolynomial in $mn$) 
total sum $N$ of entries. More precisely, we assume that $N > (mn)^7$ since the case 
of $N \leq (mn)^7$ is covered by Proposition 3.5 with a sufficiently small $\delta \leq 1/7$
(we recall that $\delta$-smooth margins are also $\delta'$-smooth with any 
$0< \delta' < \delta$).

The idea of the proof is as follows: given margins $(R,C)$ whose total sum of entries is $N$,
we construct new margins $(R', C')$ whose total sum of entries $N'$ is bounded by a polynomial
in $mn$ and a scaling map
$$\TT: \quad \Sigma(R,C) \longrightarrow \Sigma(R', C'),$$
which, roughly, scales every table $D \in \Sigma(R,C)$ by the same factor $t$.
We then deduce Theorem 1.5 for margins $(R, C)$ from that for margins $(R', C')$.

We have 
$$R' \approx t^{-1} R, \quad C' \approx t^{-1}C \quad \text{and} \quad \TT(D) \approx t^{-1} D,$$
where ``$\approx$'' stands for rounding in some consistent way.

In constructing the map $\TT$ we essentially follow the ideas of \cite{D+97}.

\subhead (4.1) Lattices, bases, and fundamental parallelepipeds \endsubhead
Let $\VV$ be a finite-dimensional real vector space and let $\Lambda \subset \VV$ 
be a {\it lattice}, that is, a discrete additive subgroup of $\VV$ which spans $\VV$. Suppose
that $\dim \VV=k$ and let $u_1, \ldots, u_k$ be a basis of $\Lambda$.
The set 
$$\Pi=\left\{\sum_{i=1}^k \lambda_i u_i: \quad 0 \leq \lambda_i <1 \quad \text{for}
\quad i=1, \ldots, k \right\}$$
is called the {\it fundamental parallelepiped} associated with the basis $u_1, \ldots, u_k$.

Suppose that $\AA$ is an affine space, with $\dim \AA=\dim \VV$, on which $\VV$ acts by 
translations: $a+v \in \AA$ for all $a \in \AA$ and $v \in \VV$ and 
$a+\left(v_1 +v_2\right)=\left(a +v _1\right)+v_2$ for all $a \in \AA$ and $v_1, v_2 \in \VV$.
Let us choose $a \in \AA$.  The set $\Lambda_a=a+\Lambda$ is called a {\it point lattice}
in $\AA$. As is known, the translations $v +\Pi: \ v \in \Lambda_a$ cover $\AA$ without overlapping.

We will also use the following standard fact. Suppose that $\Lambda_1 \supset \Lambda$ is
a finer lattice and let $|\Lambda_1/\Lambda|< \infty$ be its index. Then, for any $a, b \in \AA$ we have
$$|(a+ \Pi) \cap (b + \Lambda_1)| =|\Lambda_1/\Lambda|,$$
see for example Chapter VII of \cite{Ba02}.

Let us fix a point lattice $\Lambda_a \subset \AA$ and a fundamental parallelepiped 
$\Pi \subset \VV$ of $\Lambda$. 
 Given a point $x \in \AA$, we define its {\it rounding} 
 $y=\lfloor x\rfloor_{\Lambda_a, \Pi}$ as the unique point $y \in \Lambda_a$ such that 
$x \in  y + \Pi$.

In our case, $\VV$ is the space of real $m \times n$ matrices with the row and column 
sums equal to 0, so $\dim \VV=(m-1)(n-1)$, while $\AA$ is the affine space of $m \times n$ matrices
with prescribed integer row and column sums, so that for all $D \in \AA$ and $U \in \VV$ we have 
$D + U \in \AA$. Furthermore, let $\Lambda \subset \VV$ be the lattice of integer matrices
and let  $\Lambda'  \subset \AA$ be the point lattice consisting of integer matrices.

As is shown, for example, in \cite{D+97}, lattice $\Lambda$ has a basis consisting of the 
matrices $U_{ij}$ for $1 \leq i \leq n-1$, $1 \leq j \leq m-1$
that have $1$ in the $(i,j)$ and $(i+1, j+1)$ positions, $-1$ in the $(i+1,j)$ and $(i, j+1)$ 
positions and zeros elsewhere.  Let $\Pi$ be the fundamental parallelepiped of this basis 
$\left\{U_{ij}\right\}$. We call this parallelepiped $\Pi$ {\it standard}.
We note that
$$-2 \leq x_{ij} \leq 2 \quad \text{for all} \quad i,j \quad \text{and all} \quad X \in \Pi,
\quad X=\left(x_{ij}\right). \tag4.1.1$$

Finally, for positive integer $t$ let $\Lambda_1 = t^{-1} \Lambda$.
Hence $|\Lambda_1/\Lambda|=t^{(m-1)(n-1)}$.

\subhead (4.2) The $t$-scaling map $\TT$ \endsubhead
Let us choose a positive integer $t$ and an arbitrary $D_0 \in \Sigma(R, C)$, where
$R=\left(r_1, \ldots, r_m \right)$ and $C=\left(c_1, \ldots, c_n \right)$.
Let us define a positive $m \times n$ matrix $B$ as follows. First, we obtain $D_1$ by 
rounding up to the nearest integer every entry of $t^{-1}D_0$ and adding 2 to the result. In particular, $D_1$ is a 
positive integer matrix. Let
$$B=D_1 -t^{-1}D_0, \quad \text{so} \quad D_1=B +t^{-1}D_0.$$
 Clearly, $B=\left(b_{ij}\right)$ is an $m \times n$ matrix 
with 
$$2 \leq b_{ij} < 3 \quad \text{for all} \quad i,j. \tag4.2.1$$
Let $R'=\left(r_1', \ldots, r_m' \right)$ and $C'=\left(c_1', \ldots, c_n' \right)$ be the row and column sums of $D_1$ respectively.  Thus $R'$ and $C'$ are 
positive integer vectors and
$$\split &t^{-1} r_i + 2n \ \leq \ r_i' \ \leq \ t^{-1} r_i +3n \quad \text{for} \quad i=1, \ldots, m 
\\ &\qquad \qquad  \text{and} \\
&t^{-1} c_j + 2m \ \leq \ c_j' \ \leq \ t^{-1} c_j + 3m \quad \text{for} \quad j=1, \ldots, n. \endsplit \tag4.2.2$$
Let $\AA$ be the affine subspace of matrices with row sums $R'$ and column sums $C'$ and 
let $\Lambda' \subset \AA$ be the point lattice of integer matrices. Thus $\Lambda'=D_1 +\Lambda$,
where $\Lambda$ is the lattice of $m \times n$ integer matrices with zero row and column sums,
see Section 4.1.
For a matrix $D \in \Sigma(R, C)$ we define a matrix $\TT(D)$ by 
$$\TT(D)=\lfloor t^{-1} D + B \rfloor_{\Lambda', \Pi},$$
where $\Pi$ is the standard parallelepiped of  $\Lambda$; see Section 4.1.
In words: given a table $D \in \Sigma(R, C)$, matrix $\TT(D)$ is the unique integer matrix  such that 
the translation $\TT(D) + \Pi$ of the standard parallelepiped $\Pi$ contains $t^{-1}D + B$. 
Clearly, $\TT(D)$ is an $m \times n$ integer matrix with row sums $R'$ and column sums $C'$.
Moreover, since every entry of $t^{-1} D + B$ is at least 2 and because of (4.1.1), 
matrix $\TT(D)$ is non-negative.

Hence we have defined a map
$$\TT: \Sigma(R, C) \longrightarrow \Sigma(R', C').$$
We summarize some of its properties below.
\proclaim{(4.3) Lemma} 
\roster
\item For all $Y \in \Sigma(R', C')$ we have
$$| \TT^{-1}(Y)| \leq t^{(m-1)(n-1)};$$
\item Let $S \subset \bigl\{(i,j): \ i=1, \ldots, m,\ j=1, \ldots, n \bigr\}$ be a set of indices. Then
$$t^{-1} \sigma_S(D)  \  \leq \ \sigma_S(\TT(D))\  \leq \ t^{-1} \sigma_S(D) + 5|S|$$
for all $D \in \Sigma(R,C)$.
\endroster
\endproclaim
\demo{Proof} Given $Y \in \Sigma(R', C')$, we compute $\TT^{-1}(Y)$ as follows: we consider 
the translation $(Y-B) + \Pi$ of the standard parallelepiped $\Pi$ and observe that
$$\split \TT^{-1}(Y)=\Bigl\{ D: \quad t^{-1}D  \in (Y -B) +\Pi &
\quad \text{and} \\ &D \quad \text{is a non-negative integer matrix} \Bigr\}. \endsplit$$
Recall that $\Lambda \subset \VV$ is the lattice of $m \times n$ integer matrices
with the row and column sums equal to 0 and that $\Lambda_1=t^{-1} \Lambda$.
In the affine space of $m \times n$ matrices with row  sums $t^{-1}R$ and column sums $t^{-1}C$ let
us consider the point lattice $\Lambda_1'=t^{-1}D_0 + \Lambda_1$ consisting of matrices $t^{-1}D$
where $D$ is an integer matrix. Then 
$$|((Y-B) + \Pi) \cap \Lambda_1'| =|\Lambda_1/\Lambda|=t^{(m-1)(n-1)}$$
and Part (1) follows. 
Part (2) follows because of (4.1.1) and (4.2.1).
{\hfill \hfill \hfill} \qed
\enddemo

\proclaim{(4.4) Lemma} Suppose that 
$$r_i', \  c_j' \ \geq \ (mn)^2 \quad \text{for all} \quad i,j.$$
Then, for any $\zeta \geq 0$ we have 
$$\split &\Pr \Bigl\{D \in \Sigma(R,C): \ \sigma_S(D) \geq t\zeta \Bigr\} \ \leq \ \beta 
\Pr \Bigl\{Y \in \Sigma(R',C'): \ \sigma_S(Y) \geq \zeta \Bigr\} \\ 
&\qquad \qquad \text{and} \\
&\Pr \Bigl\{D \in \Sigma(R,C): \ \sigma_S(D) \leq t\zeta \Bigr\} \ \leq \ \beta 
\Pr \Bigl\{Y \in \Sigma(R',C'): \  \sigma_S(Y) \leq \zeta +5|S| \Bigr\},\endsplit$$
where $\beta > 0$ is an absolute constant.
\endproclaim
\demo{Proof} By Part (2) of Lemma 4.3, if $\sigma_S(D) \geq t \zeta$ then 
$\sigma_S(Y) \geq \zeta$ for $Y=\TT(D)$. Using Part (1) of Lemma 4.3, we can write
$$\split \Pr \Bigl\{D \in \Sigma(R,C):\ 
 &\sigma_S(D)  \geq t \zeta \Bigr\} = {\big|D \in \Sigma(R, C): \ 
\sigma_S(D) \geq t \zeta \big| \over  |\Sigma(R,C)|} \\ \leq \ & t^{(m-1)(n-1)}
 { \big|Y \in \Sigma(R', C'): \ \sigma_S(Y) \geq \zeta \big|\over |\Sigma(R,C)|} \\
=&{|\Sigma(R', C')| \over |\Sigma(R, C)|} t^{(m-1)(n-1)}   \Pr \Bigl\{Y \in \Sigma(R', C'): \ 
\sigma_S(Y) \geq \zeta \Bigr\}. \endsplit $$ 
Similarly, by Part (2) of Lemma 4.3, if $\sigma_S(D) \leq t \zeta$ then $\sigma_S(Y) \leq \zeta + 5|S|$
for $Y=\TT(D)$ and 
$$ \split &\Pr \Bigl\{D \in \Sigma(R,C):\ \sigma_S(D)  \leq t \zeta \Bigr\} \\ &\qquad \qquad \leq \ 
 {|\Sigma(R', C')| \over |\Sigma(R, C)|} t^{(m-1)(n-1)}   \Pr \Bigl\{Y \in \Sigma(R', C'): \ 
\sigma_S(Y) \leq \zeta +5|S| \Bigr\}. \endsplit$$
It is shown in \cite{D+97} that for sufficiently large margins, the number of contingency 
tables is approximated within a constant factor by the volume of the corresponding transportation polytope; see Section 1.2. In particular, estimates of \cite{D+97} imply that 
$$|\Sigma(R', C')| \ \leq \ \beta_1 \vl \PP(R', C')  \quad \text{and} \quad 
|\Sigma(R, C)| \ \geq \ \beta_2 \vl \PP(R, C)$$
for some absolute constants $\beta_1, \beta_2  >0$.

From (4.2.2), we have 
$$\split &r_i \ \geq \ t(r_i' - 3n) \ \geq \ t r_i' \left(1 -{3 \over m^2 n} \right) \quad \text{for} \quad i=1, \ldots, m\quad \text{and} \\
&c_j \ \geq \ t(r_i' - 3m) \ \geq \ t c_j' \left(1 -{3 \over m n^2} \right) \quad \text{for} \quad j=1, \ldots, n.
\endsplit$$
It follows then that 
$$\vl \PP(R,C) \ \geq \ \beta_3 t^{(m-1)(n-1)} \vl \PP(R', C')$$
for some absolute constant $\beta_3>0$.
The result now follows.
{\hfill \hfill \hfill} \qed
\enddemo 

Next, we show that the $t$-scaling map $\TT$ almost scales the typical table provided 
the margins $R', C'$ are large enough, that is, $Z' \approx t^{-1} Z$. The idea of the proof is roughly the following: 
if margins $(R', C')$ and $(R, C)$ are large enough, then the corresponding typical 
tables $Z'$ and $Z$ roughly optimize the functional $\sum_{i,j} \ln x_{ij}$ on the corresponding 
transportation polytopes and hence the map $X \longmapsto t X$ roughly maps $Z'$ to $Z$.

\proclaim{(4.5) Lemma} Let $Z=\left(z_{ij}\right)$ be the typical table with margins $(R,C)$, let 
$Z'=\left(z_{ij}' \right)$ be the typical table with margins $(R', C')$ obtained by $t$-scaling and suppose that 
$$z_{ij}' \ \geq \ (mn)^4+3 \quad \text{for all} \quad i,j.$$
Then 
$$ \Big|{z_{ij} \over t z_{ij}'} -1 \Big| \ \leq \ {\beta \over mn} \quad \text{for all} \quad i,j$$
and some absolute constant $\beta>0$.
\endproclaim

\demo{Proof} First, we prove some useful inequalities for the function 
$$g(x)=(x+1) \ln (x+1) - x \ln x.$$
We have 
$$g(tx) -g(x) =\int_x^{tx} g'(y) \ dy = \int_x^{tx} \ln \left({y+1 \over y}\right) \ dy 
\leq \int_x^{tx} {dy \over y} =\ln (tx) -\ln x=\ln t.$$
Also,
$$\split g(x) =&(x+1) \ln (x+1) - (x+1) \ln x +(x+1)\ln x - x\ln x \\=&(x+1) \ln \left({x+1 \over x}\right) + \ln x=
\ln x + 1 + O\left({1 \over x}\right) \quad \text{for} \quad x \geq 1. \endsplit $$
Finally, we note that 
$$g''(x)=-{1 \over x(x+1)}.$$

Since from (4.2.2) we have 
$$r_i \ \leq \ t r_i' \quad \text{and} \quad c_j \ \leq \ t c_j' \quad \text{for all} \quad i,j$$
we have 
$$\max \Sb X \in \PP(R, C) \endSb g(X)  \ \leq \ \ln t + \max \Sb X \in \PP(R', C') \endSb g(X). \tag4.5.1$$
Let $B$ be the matrix constructed in Section 4.2 and let $W=t(Z' -B) \in \PP(R, C)$. Hence
$$w_{ij} \ \geq \ t (mn)^4 \quad \text{for all} \quad i,j.$$
Since
$$g\left(w_{ij}\right) =1 + \ln w_{ij} + O\left({1 \over m^4 n^4}\right) \quad \text{and} \quad
g\left(z_{ij}'\right) =1 + \ln z_{ij}' + O \left({1 \over m^4 n^4} \right),$$
we have 
$$g(W) =g(Z') + \ln t + O\left({1 \over m^3n^3}\right).$$
From (4.5.1) it follows that
$$g(Z)-g(W)=O\left({1 \over m^3n^3} \right). \tag4.5.2$$
Next, we are going to exploit the strong concavity of $g$ and use the following standard inequality: 
\smallskip
if $g''(x) \leq -\alpha$ for some $\alpha>0$ and all $a \leq x \leq b$ then
$$g\left({a+b \over 2}\right)-{1 \over 2} g(a) -{1 \over 2} g(b) \ \geq \  {\alpha (b-a)^2 \over 8}.$$
\smallskip
If for some $i,j$ we have $|w_{ij}-z_{ij}| \geq (mn)^{-1} w_{ij}$, then in view of (4.5.2), for some point $U$ on the interval 
connecting $W$ and $Z$ and all sufficiently large $mn$, we will have 
$$g(U) > g(Z),$$
which is a contradiction.
Thus 
$$\left|{z_i \over w_{ij}}-1 \right| \ \leq \ {1 \over mn} \quad \text{for all} \quad i,j$$
and all sufficiently large $mn$.
Since
$$\left|{w_{ij} \over tz_{ij}' }-1 \right|  \ \leq \ {3 \over z_{ij}'} \ \leq \ {3 \over (mn)^4},$$
 the proof follows.
{\hfill \hfill \hfill} \qed
\enddemo

\subhead (4.6) Proof of Theorem 1.5 \endsubhead
Without loss of generality we assume that $N \geq (mn)^7$ since the case of a polynomially bounded 
$N$ is handled in Proposition 3.5.

Let us choose 
$$t=\left\lfloor {N \over (mn)^6} \right\rfloor$$
and consider the $t$-scaling map 
$\TT: \Sigma(R, C) \longrightarrow \Sigma(R', C')$ .
Since margins $(R,C)$ are $\delta$-smooth, we have 
$$(mn)^6 \ \leq \ N' \ \leq \ (mn)^7 \quad \text{and} \quad 
r_i', c_j' \ \geq \ (mn)^4 \quad \text{for all} \quad i,j$$
and all sufficiently large $n \geq m$.

Let us choose $0 < \delta_1 < \delta$. It follows by 
(4.2.2) that the margins $(R', C')$ are $\delta_1$-smooth
for all sufficiently large $n \geq m$.
Let $Z'$ be the typical table of $(R', C')$, $Z'=\left(z_{ij}'\right)$.
By Corollary 2.5,
$$z_{ij}'  \ \geq \ (\delta_1)^3 {N' \over mn}.$$
Therefore, for all sufficiently large $m+n$ we have
$$z_{ij}' \ \geq \ (mn)^4 + 3.$$
The result now follows by Lemmas 4.4, 4.5, and Proposition 3.5 applied to $(R', C')$.
 {\hfill \hfill \hfill} \qed
 
\head Acknowledgment \endhead

The author is grateful to Alexander Yong for asking in what sense 
the ``typical table'' introduced in \cite{Ba09} and named so in \cite{B+08} was typical. This paper is an attempt to answer that question. After the first version of this paper was written, John Hartigan 
pointed out to connections with the maximum entropy principle and
suggested Theorem 1.7 to the author (cf. \cite{BH09}), 
which led to a substantial simplification of the original proof and some strengthening of the main result, Theorem 1.5.


\Refs
\widestnumber\key{AAAA}

\ref\key{Ba02}
\by A. Barvinok
\book A Course in Convexity
\bookinfo Graduate Studies in Mathematics, 54
\publ American Mathematical Society
\publaddr Providence, RI
\yr 2002
\endref

\ref\key{Ba08}
\by A. Barvinok
\paper On the number of matrices and a random matrix with prescribed row and column sums and 0-1 entries
\paperinfo preprint arXiv:0806.1480 
\yr 2008
\endref
	
\ref\key{B+08}
\by A. Barvinok, Z. Luria, A. Samorodnitsky, and A. Yong
\paper An approximation algorithm for counting contingency tables
\paperinfo preprint {\tt arXiv:0803.3948}
\jour Random Structures $\&$ Algorithms, to appear 
\yr 2008
\endref

\ref\key{Ba09}
\by A. Barvinok
\paper Asymptotic estimates for the number of contingency tables, integer flows, and volumes of transportation polytopes
\jour International Mathematics Research Notices
\vol 2009
\yr 2009
\pages 348--385
\endref

\ref\key{BH09}
\by A. Barvinok and J.A. Hartigan
\paper Maximum entropy Gaussian approximation for the number of integer points and volumes of polytopes
\paperinfo preprint {\tt arXiv:0903.5223}
\yr 2009
\endref

\ref\key{C+06}
\by M. Cryan, M. Dyer, L.A. Goldberg,  M. Jerrum, and R. Martin
\paper Rapidly mixing Markov chains for sampling contingency tables with a constant number of rows
\jour SIAM Journal on Computing 
\vol 36 
\yr 2006
\pages 247--278
\endref


\ref\key{DE85}
\by P. Diaconis and B. Efron
\paper Testing for independence in a two-way table: new interpretations of the chi-square statistic. With discussions and with a reply by the authors
\jour The Annals of Statistics 
\vol 13 
\yr 1985
\pages 845--913
\endref

\ref\key{DG95}
\by P. Diaconis and A. Gangolli
\paper Rectangular arrays with fixed margins
\inbook Discrete probability and algorithms (Minneapolis, MN, 1993)
\pages 15--41 
\bookinfo The IMA Volumes in Mathematics and its Applications, 72
\publ  Springer
\publaddr New York
\yr 1995
\endref

\ref\key{D+97} 
\by M. Dyer, R. Kannan, and J. Mount 
\paper Sampling contingency tables 
\jour Random Structures $\&$ Algorithms
\vol 10 
\yr 1997 
\pages 487--506 
\endref

\ref\key{Go63}
\by I. J. Good
\paper Maximum entropy for hypothesis formulation, especially for multidimensional contingency tables
\jour The Annals of Mathematical Statistics
\vol 34
\yr 1963
\pages 911-934  
\endref

\ref\key{GM08}
\by C. Greenhill and B.D. McKay
\paper Asymptotic enumeration of sparse nonnegative integer matrices with specified row and column sums
\jour  Advances in Applied Mathematics
\vol 41 
\yr 2008
\pages 459--481
\endref

\ref\key{Le01}
\by M. Ledoux
\book The Concentration of Measure Phenomenon
\bookinfo Mathematical Surveys and Monographs, 89
\publ American Mathematical Society
\publaddr Providence, RI
\yr 2001
\endref

\ref\key{Ne69}
\by P.E. O'Neil
\paper Asymptotics and random matrices with row-sum and column-sum restrictions
\jour  Bulletin of the American Mathematical Society 
\vol 75 
\yr 1969 
\pages 1276--1282
\endref

\ref\key{NN94}
\by Yu. Nesterov and A. Nemirovskii
\book Interior-Point Polynomial Algorithms in Convex Programming
\bookinfo SIAM Studies in Applied Mathematics, 13
\publ Society for Industrial and Applied Mathematics (SIAM)
\publaddr Philadelphia, PA
\yr 1994
\endref


\endRefs
\enddocument
\end